\documentclass[12pt]{article}
\usepackage[final]{epsfig}
\usepackage{graphics}
\usepackage{amsmath}
\usepackage{amsfonts}
\usepackage{latexsym}
\usepackage{amssymb}
\usepackage{graphicx}
\newtheorem{lemma}{Lemma}
\newtheorem{proposition}[lemma]{Proposition}
\newtheorem{remark}[lemma]{Remark}

\newtheorem{theorem}[lemma]{Theorem}

\begin{document}
\newcommand{\eps}{{\varepsilon}}
\newcommand{\proofend}{\hfill$\Box$\bigskip}
\newcommand{\C}{{\mathbb C}}
\newcommand{\Q}{{\mathbb Q}}
\newcommand{\R}{{\mathbb R}}
\newcommand{\T}{{\mathbb T}}
\newcommand{\Z}{{\mathbb Z}}
\newcommand{\RP}{{\mathbf {RP}}}
\def\proof{\paragraph{Proof.}}

%%%%%%%%%%%%%%%%%%%%%%%%%%%%%%%%%%%%%%%%%%%%%%%%%%%%%%%%%%

\title {Birkhoff billiards are insecure}
\author{Serge Tabachnikov\thanks{
Department of Mathematics,
Pennsylvania State University, University Park, PA 16802, USA;
e-mail: \tt{tabachni@math.psu.edu}
}
\\
}
\date{\today}
\maketitle
\begin{abstract}
We prove that every compact plane billiard, bounded by a smooth curve, is insecure: there exist pairs of points $A,B$ such that no finite set of points can block all billiard trajectories from $A$ to $B$.
\end{abstract}

Two points $A$ and $B$ of a Riemannian manifold $M$ are called {\it secure} if there exists a finite set of points $S\subset M-\{A,B\}$ such that every geodesic connecting $A$ and $B$ passes through a point of $S$. One says that the set $S$ blocks $A$ from $B$. 
A manifold is called secure (or has the {\it finite blocking property}) if any pair of its points is secure.  For example, every pair of non-antipodal points of  the Euclidean sphere is secure, but a pair of antipodal points is not secure, so the sphere is insecure. A flat torus of any dimension is secure.

In the recent years, the notion of security has attracted a considerable attention, see \cite{B-G,Gut,G-S,L-S,Mo1,Mo2,Mo3,S-S}. This notion extends naturally to Riemannian manifolds with boundary, in which case one considers billiard trajectories from $A$ to $B$ with the billiard reflection off the boundary. 

In this note we consider a compact plane billiard domain $M$ bounded by a smooth curve  and prove that $M$ is insecure. More specifically, one has the following local insecurity result. Consider a sufficiently short outward convex arc $\gamma \subset \partial M$ with end-points $A$ and $B$ (such an arc always exists).

\begin{theorem} \label{main}
The pair $(A,B)$ is insecure.
\end{theorem}

\proof
Denote by $T_n$ the polygonal line $A=P_0,P_1,\dots P_{n-1},P_n=B,\ \ P_i\in\gamma$,  of minimal length; this is a billiard trajectory from $A$ to $B$. If $n$ is large then $T_n$ lies in a small neighborhood of $\gamma$.

Working toward contradiction, assume that a finite set of points $S\subset M-\{A,B\}$ blocks every billiard trajectory from $A$ to $B$. Decompose $S$ as $S'\cup S''$ where the points of $S'$ lie on the boundary and the points of $S''$ lie inside the billiard table. For $n$ large enough, the trajectory $T_n$ is disjoint from $S''$. We want to show that there is a sufficiently large $n$ such that the set ${\cal P}_n=\{P_1,\dots P_{n-1}\}$ is disjoint from $S'$.

Let $s$ be the arc-length parameter and $k(s)$ the curvature of $\gamma$. Let $\sigma$ be a new parameter on the arc $\gamma$ such that $d \sigma = (1/2)k^{2/3} ds$. 
By rescaling the arc $\gamma$, we may assume that the range of $\sigma$ is $[0,1]$ with $\sigma(0)=A$ and $\sigma(1)=B$.  Let $Q_0=A,Q_1,\dots,Q_{n-1},Q_n=B$ be the points that divide the $\sigma$-measure of $\gamma$ into $n$ equal parts, that is, $Q_m=\sigma(m/n)$.

\begin{proposition} \label{distr}
One has: $|P_n - Q_n| = O(1/n^2)$.
\end{proposition}

\begin{remark}
{\rm This Claim is consistent with Theorem 6 (iii) of \cite{M-V} which describes the limit distribution of the vertices of the inscribed polygons that best approximate a convex curve relative the deviation of the perimeter length. 
}
\end{remark}

To prove Proposition \ref{distr}, we use the theory of interpolating Hamiltonians, see \cite{Me1,Me2} and especially \cite{M-M}. Recall the relevant facts from this theory.

First, some generalities about plane billiards (see, e. g., \cite{Ta1,Ta2}).
The phase space $X$ of the billiard ball map consists of inward unit tangent vectors $(x,v)$ to $M$ with the foot point $x$ on the boundary  $\partial M$;  $x$ is the position of the billiard ball and $v$ is its velocity. The billiard ball map $F$ takes $(x,v)$ to the vector obtained by moving $x$ along $v$ until it hits $\partial M$ and then elastically reflecting $v$ according to the law ``angle of incidence equals angle of reflection".  Let $\phi$ be the angle made by $v$ with the positive direction of $\partial M$. Then $(s,\phi)$ are coordinates in $X$. The area form $\omega=\sin\phi\ d\phi\wedge ds$ is $F$-invariant.

In a nutshell, the theory of interpolating Hamiltonians asserts that the billiard ball map equals an integrable symplectic map, modulo smooth symplectic maps that fix the boundary of the phase space $X$ to all orders. More specifically,  one can choose new symplectic coordinates $H$ and $Z$ near the boundary $\phi=0$ such that $\omega=dH \wedge dZ$, $H$ is an integral of the map $F$, up to all orders in $\phi$, and 
\begin{equation} \label{shift}
F^*(Z)=Z+H^{1/2},
\end{equation}
also up to all orders in $\phi$. The function $H$ is given by a series in even powers of $\phi$, namely, 
\begin{equation} \label{Ham}
H=k^{-2/3} \phi^2 + O(\phi^4),
\end{equation}
 and this series is uniquely determined by the above conditions on $H$ and $Z$. 

\begin{lemma} \label{choice}
One may choose the coordinate $Z$ in such a way  that
$Z=\sigma + O(\phi^2)$.
\end{lemma}

\proof
Let $Z=f(s)+g(s)\phi+O(\phi^2)$. We have: $\omega=dH \wedge dZ$. Equating the coefficients of $\phi\ d\phi\wedge ds$ and of $\phi^2\ d\phi\wedge ds$ and using (\ref{Ham}) we obtain the equations:
$$
2k^{-2/3}(s) f'(s)=1,\quad 2k^{-2/3}(s) g'(s)+\frac{2}{3}k^{-5/3}(s)k'(s) g(s)=0.
$$
The first equation implies that $df=d\sigma$ and the second that $g=C k^{-1/3}$ where $C$ is a constant. We can choose $f(0)=0$. Since $Z$ is defined up to summation with functions of $H$, it follows from (\ref{Ham}) that the term $g(s)\phi$ can be eliminated by subtracting $C H^{1/2}$.
\proofend

Now we can prove Proposition \ref{distr}. The billiard trajectory $T_n$ corresponds to a phase orbit $x_0,\dots,x_n,\ F(x_i)=x_{i+1}$. Since $H$ is an integral of the map $F$, the orbit $x_0,\dots,x_n$ lies on a level curve $H=c_n$. Due to (\ref{shift}), we have: $n \sqrt{c_n} = O(1)$, and hence $c_n=O(1/n^2)$ which, in view of (\ref{Ham}), implies that 
\begin{equation} \label{1/n}
\phi = O\left(\frac{1}{n}\right). 
\end{equation}
Consider $\sigma$ and $Z$ as  functions on the phase space $X$. Since $\sigma(x_m)=P_m$ and the $\sigma$-coordinate of $Q_m$ is $m/n$, we need to show that  
\begin{equation} \label{rat}
\sigma (x_m)=\frac{m}{n}+O\left(\frac{1}{n^2}\right).
\end{equation}
Since $F$ is a shift in $Z$-coordinate, see (\ref{shift}), one has:
$$ 
Z(x_m)=\frac{m}{n} \left( Z(x_n)-Z(x_0) \right)= \frac{m}{n} \left( \sigma(x_n)-\sigma(x_0) \right)+O\left(\frac{1}{n^2}\right)=\frac{m}{n}+O\left(\frac{1}{n^2}\right),
$$
the second equality due to Lemma \ref{choice} and (\ref{1/n}). 
This proves Proposition \ref{distr}.
 \proofend
 
From now on, we identify the arc $\gamma$ with the segment $[0,1]$ using the parameter $\sigma$; the points $P_1,\dots,P_{n-1}$ are considered as reals between $0$ and $1$. Assume that a finite set $S'=\{t_1,\dots,t_k\}\subset (0,1)$ is blocking, that is, for all sufficiently large  $n$, one has ${\cal P}_n \cap S' \neq \emptyset$. 

Some of the numbers $t_i\in S'$ may be rational; denote them by   $p_i/q_i,\ i=1,\dots,l$ (fractions in lowest terms), and let $Q=q_1\cdots q_l$. Set $n_i=1+(N+i)Q,\ i=0,\dots,k$. 

\begin{proposition} \label{void}
For $N$ sufficiently large, at least one of the sets ${\cal P}_{n_i}$ is disjoint from $S'$.
\end{proposition}

\proof
Assume not. Then, by the Pigeonhole Principle, there exist $l,i,j$ such that 
$t_l\in {\cal P}_{n_i} \cap {\cal P}_{n_j}.$  According to Proposition \ref{distr}, there is a constant $C$ (independent of $n$) such that, for $P_m \in {\cal P}_n$, one has: 
$$
\left|P_m - \frac{m}{n}\right|<\frac{C}{n^2}. 
$$
Therefore 
\begin{equation} \label{ineq}
\left|t_l - \frac{m_1}{n_i}\right|<\frac{C}{n_i^2},\ \ \left|t_l - \frac{m_2}{n_j}\right|<\frac{C}{n_j^2}
\end{equation} 
for some $m_1, m_2$.

\begin{lemma} \label{irrt}
If $N$ sufficiently large then $t_i \notin \Q$.
\end{lemma}

\proof
First, we claim that, given  a fraction $p/q$ and a constant $C$, if 
$$
\left|\frac{p}{q}-\frac{m}{n}\right|<\frac{C}{n^2}
$$
for all sufficiently large $n$ then $m/n= p/q$.

Indeed, if $m/n\neq p/q$ then $1\leq |pn-qm|$, hence 
$$
\frac{1}{qn}\leq \left|\frac{p}{q}-\frac{m}{n}\right|<\frac{C}{n^2},
$$
which cannot hold for $n>Cq$.

Next, we claim that, for all  $M,N\in \Z$ and each $i=1,\dots,l$, 
$$
\frac{M}{1+NQ}\neq \frac{p_i}{q_i}.
$$
Indeed, if the equality holds then $Mq_i=p_i(1+NQ)$. The right hand side is divisible by $q_i$ but $1+NQ$ is coprime with $q_i$; this contradicts the assumption that $q_i$ and $p_i$ are coprime.

The two claims combined imply the lemma.
\proofend

Next,  (\ref{ineq}) and the triangle inequality imply that 
$$
\left|\frac{m_1}{n_i} - \frac{m_2}{n_j}\right| < C\left( \frac{1}{n_i^2}+ \frac{1}{n_j^2}\right)
$$
for some $m_1,m_2$. It follows that 
$$
|m_1n_j-m_2n_i|<C\left( \frac{n_j}{n_i}+ \frac{n_i}{n_j}\right).
$$
The expression in the parentheses on the right hand side has limit 2, as $N\to\infty$, 
hence one has, for sufficiently great $N$,
\begin{equation} \label{ineq2}
|m_1n_j-m_2n_i|<3C.
\end{equation}

Denote by ${\cal M}$ the (finite) set of fractions with the denominators $jQ,\ j\in \{1,2,\dots,k\}$, and let $\delta >0$ be the distance between the sets $S'-\Q$ and ${\cal M}$.

\begin{lemma} \label{notsmall}
For sufficiently large $N$, one has: 
$$
|m_1n_j-m_2n_i|> \delta Q^2N/2.
$$
\end{lemma}
 
\proof
For $N$ large enough, it follows from (\ref{ineq}) that 
$$
\left|t_l - \frac{m_1}{n_i}\right|<\frac{\delta}{2}.
$$
Since $t_l \notin \Q$, it follows that the distance from $m_1/n_i$ to ${\cal M}$ is greater than $\delta/2$. One has:
$$
|m_1n_j-m_2n_i|=|n_j-n_i|\ n_i \left| \frac{m_1}{n_i}-\frac{m_2-m_1}{n_j-n_i}\right| > Q\cdot QN\cdot \frac{\delta}{2},
$$
as claimed.
\proofend

Finally,  for large $N$, Lemma \ref{notsmall} contradicts inequality (\ref{ineq2}), and Proposition \ref{void} follows.
\proofend

This Proposition implies Theorem \ref{main}, and we are done. 
\proofend

\begin{remark}
{\rm Theorem \ref{main}, along with its proof, can be extended to  billiards in higher dimensional Euclidean spaces: the role of the curve $\gamma$ is played by the shortest geodesic on the boundary of the billiard table connecting $A$ and $B$.
}
\end{remark}

\bigskip

{\bf Acknowledgments}. Many thanks to K. Burns, R. Schwartz and L. Stojanov for their interest. The  author was partially supported by an NSF grant DMS-0555803.


\bigskip

\begin{thebibliography}{99}

\bibitem{B-G} K. Burns, E. Gutkin. {\it Growth of the number of geodesics between points and insecurity for riemannian manifolds}. Preprint arXiv:math/0701579

\bibitem{Gut} E. Gutkin. {\it Blocking of billiard orbits and security for polygons and flat surfaces}. Geom. Funct. Anal. {\bf 15} (2005), 83--105.

\bibitem{G-S} E. Gutkin, V. Schroeder. {\it Connecting geodesics and security of configurations in compact locally symmetric spaces}. Geom. Dedicata {\bf 118} (2006), 185--208. 

\bibitem{L-S} J.-F. Lafont, B. Schmidt. {\it Blocking light in compact Riemannian manifolds}. Geometry and Topology, to appear,  arXiv:math/0607789

\bibitem{M-M} S. Marvizi, R. Melrose. {\it Spectral invariants of convex planar regions. }J. Diff. Geom. {\bf 17} (1982), 475-502.

\bibitem{M-V} D. McClure, R. Vitale.
{\it Polygonal approximation of plane convex bodies.} 
J. Math. Anal. Appl. {\bf 51} (1975),  326--358. 

\bibitem{Me1} R. Melrose. {\it Equivalence of glancing hypersurfaces.} Invent. Math. {\bf 37} (1976), 165-192. 
 
\bibitem{Me2} R. Melrose. {\it Equivalence of Glancing Hypersurfaces 2.} Math. Ann. {\bf 255} (1981), 159-198.

\bibitem{Mo1} T. Monteil. {\it A counter-example to the theorem of Hiemer and Snurnikov}. J. Statist. Phys. {\bf 114} (2004),  1619--1623.

\bibitem{Mo2} T. Monteil. {\it On the finite blocking property}. Ann. Inst. Fourier  {\bf 55} (2005),  1195--1217.

\bibitem{Mo3} T. Monteil. {\it Finite blocking property versus pure periodicity}. Preprint
arXiv:math/0406506

\bibitem{S-S} B. Schmidt, J. Souto. {\it Chords, light, and another synthetic characterization of the round sphere}. Preprint arXiv:0704.3642

\bibitem{Ta1} S. Tabachnikov. {\it Billiards}. Soc. Math. de France,
Paris, 1995.

\bibitem{Ta2} S. Tabachnikov. {\it Geometry and billiards}.  Amer. Math. Soc., Providence, RI, 2005.

\end{thebibliography}
\end{document}